\documentclass[10pt]{book}
\usepackage{array}
\usepackage{setspace}
\usepackage[sort&compress]{natbib}
\bibpunct{(}{)}{;}{a}{}{,}
\usepackage{graphicx}
\usepackage{caption,subcaption,fancyhdr}
\usepackage{amssymb}
\usepackage{amsmath}
\usepackage{epstopdf}
\usepackage{amsthm}
\usepackage{enumerate}

\usepackage{breakcites}
\usepackage{hyperref}

\makeatletter
\newcommand*{\rom}[1]{\expandafter\@slowromancap\romannumeral #1@}
\makeatother

\newtheorem{theorem}{Theorem}

\newtheorem{definition}{Definition}
\newtheorem{proposition}{Proposition}

\newcommand{\omegaopt}{\omega_{\textrm{opt}}}
\newcommand{\Ropt}{R_{\textrm{opt}}}
\newcommand{\omegaBH}{\omega_{\textrm{BH}}}
\newcommand{\omegaGW}{\omega_{\textrm{GW}}}
\newcommand{\omegahat}{\hat{\omega}}

\begin{document}
\pagestyle{fancy}
\renewcommand{\baselinestretch}{1.2}
\lhead[\fancyplain{} \leftmark]{}
\chead[]{}
\rhead[]{\fancyplain{}\rightmark}
\cfoot{}

\markright{
\hbox{\footnotesize\rm Statistica Sinica
}\hfill\\[-13pt]
\hbox{\footnotesize\rm
}\hfill }

\markboth{\hfill{\footnotesize\rm XUEYING TANG, KE LI AND MALAY GHOSH}\hfill}
{\hfill {\footnotesize\rm ABOS OF POLYNOMIAL-TAILED DISTRIBUTIONS} \hfill}
\renewcommand{\thefootnote}{}
$\ $\par \fontsize{10.95}{14pt plus.8pt minus .6pt}\selectfont
\vspace{0.8pc}
\centerline{\large\bf BAYESIAN MULTIPLE TESTING UNDER SPARSITY}
\vspace{2pt}
\centerline{\large\bf FOR POLYNOMIAL-TAILED DISTRIBUTIONS}
\vspace{.4cm}
\centerline{Xueying Tang, Ke Li and Malay Ghosh}
\vspace{.4cm}
\centerline{\it University of Florida, Southwestern University of Finance and Economics}
\centerline{\it and University of Florida}
\vspace{.55cm}
\fontsize{9}{11.5pt plus.8pt minus .6pt}\selectfont

\doublespacing
\begin{quotation}
\noindent {\it Abstract:}
This paper considers Bayesian multiple testing under sparsity for polynomial-tailed distributions satisfying a monotone likelihood ratio property. Included in this class of distributions are the Student's t, the Pareto, and many other distributions. We prove some general asymptotic optimality results under fixed and random thresholding. As examples of these general results, we establish the Bayesian asymptotic optimality of several multiple testing procedures in the literature for appropriately chosen false discovery rate levels. We also show by simulation that the Benjamini-Hochberg procedure with a false discovery rate level different from the asymptotically optimal one can lead to high Bayes risk.
\par

\vspace{9pt}
\par
\end{quotation}\par

\fontsize{10.95}{14pt plus.8pt minus .6pt}\selectfont
\setcounter{chapter}{1}
\setcounter{equation}{0} 
\noindent {\bf 1. Introduction}
\par \smallskip

Multiple testing has become a topic of growing importance in recent years. Its importance is particularly felt in the event of inference under sparsity, detecting a few signals in the midst of multiple noises. Applications abound, for example, in genetics, engineering, biology, and finance, just to name a few. A specific example is when one needs to identify a handful of genes attributable to a certain disease in the midst of thousands of others.

Currently, the most widely used approach for multiple testing is the one due to \citet*{benjamini1995controlling} that controls the false discovery rate (FDR). Since then, there are several noteworthy contributions in the area. Among others, we refer to \citet*{benjamini2001control}, \citet*{efron2002empirical}, \citet*{abramovich2006special}, 
\citet*{donoho2006asymptotic}, \citet*{gavrilov2009adaptive}, \citet*{genovese2002operating,genovese2004stochastic}, \citet*{sarkar2002some}, and \citet*{storey2002direct}. 

Recently, \citet*{bogdan2008comparison} conducted an extensive simulation study to find closeness of the Benjamini-Hochberg procedure to an optimal Bayes procedure for multiple hypothesis testing under normality of the data. Later, in \citet*{bogdan2011asymptotic} and \citet*{frommlet2013some}, it was shown that several multiple testing procedures, including the Benjamini-Hochberg procedure, asymptotically attained the Bayes oracle property under sparsity, once again under normality. As an extension of \citet*{bogdan2011asymptotic}, \citet*{neuvial2012false} studied properties of FDR thresholding with observations coming from the Subbotin family that includes Laplace and normal distributions as special cases.

Here we consider distributions with polynomial tails. As we show later, the Bayes rules of the multiple testing problem of normal distributions and polynomial-tailed distributions produce quite different Bayes risks. Both type \rom{1} and type \rom{2} errors play a role in the limiting (as the number of tests goes to infinity) Bayes risk of the oracle multiple testing procedure for polynomial-tailed distributions while, as shown in \citet*{bogdan2011asymptotic}, the Bayes risk is asymptotically determined solely by type \rom{2} errors for normal distributions. As indicated by \citet*{chi2007performance}, in multiple testing problems for polynomial-tailed distributions, controlling FDR under a certain threshold level $\alpha^*$ leads to asymptotically zero power. As a result, a vanishing FDR level is very unlikely to define an asymptotically optimal procedure in terms of Bayes risk. As this is not the case for normal distributions, we were motivated to study the performance of the Benjamini-Hochberg procedure and some other multiple testing procedures for polynomial-tailed distributions. We study the asymptotic optimality of multiple testing procedures for a general class of such distributions, including Student's t, Pareto, and many others.

Our framework follows that of \citet*{bogdan2011asymptotic} and \citet*{neuvial2012false}, where the multiple testing problem is addressed in a decision theoretic framework. Suppose there are $m$ independent observations, each of which comes from a mixture of two distributions in the scale family of a polynomial-tailed distribution. We are interested in testing simultaneously which distribution each observation comes from. We assume the loss of wrong decisions for the $m$ tests is the sum of the losses of wrong decisions for individual tests (\citet*{lehmann1957theorya,lehmann1957theoryb}). For each test, nonzero loss occurs if and only if a type \rom{1} or a type \rom{2} error is made. All our results are obtained in an asymptotic framework that ensures that the limiting power of an individual test based on the Bayes oracle threshold converges to a constant between zero and one. 

After finding the oracle Bayes risk, we define asymptotic Bayesian optimality under sparsity (ABOS) analogous to \citet*{bogdan2011asymptotic}. Under the asymptotic framework, a necessary and sufficient condition is provided for a fixed thresholding procedure to be ABOS. A single constraint guarantees that the risks from the two error types converge to the optimal risk for our proposed class of distributions. 

A more practically meaningful result is that we obtain a general sufficient condition for a random thresholding procedure to obtain ABOS. The condition requires comparison of a random threshold with a fixed ABOS threshold that is sometimes easier to work with than bounding the Bayes risk directly, as in \citet*{bogdan2011asymptotic}. Our general results show that the procedures controlling Bayesian false discovery rate (BFDR), the procedure of \citet*{genovese2002operating}, and the Benjamini-Hochberg procedure are all ABOS if the FDR level is chosen properly. On the other hand, it is shown via simulation that the Benjamini-Hochberg procedure with an FDR level different from the optimal one can lead to high Bayes risk.


The remaining sections of this article are organized as follows. In Section 2, we describe our asymptotic framework. The Bayes oracle rule and its Bayes risk are derived in this section, and we list a few important distributions with polynomial tails. Two general results about the conditions for fixed or random thresholding procedures to be ABOS are presented in Section 3. In Section 4, we provide conditions under which several procedures, including Benjamini-Hochberg, are ABOS by applying the results in Section 3. Section 5 contains numerical results suggesting the non-optimality of the Benjamini-Hochberg procedure. Some final remarks are made in Section 6. Proofs of theoretical results are provided in the supplementary material.
\par \bigskip

\setcounter{chapter}{2}
\setcounter{equation}{0} 
\noindent {\bf 2. Oracle Bayes Risk and Asymptotic Framework}
\par \smallskip
Suppose we have $m$ independent observations $\mathbf{X} = (X_1, \ldots, X_m)$ from the same distribution. Let $D$ and $d$ be the cumulative distribution function (cdf) and the probability density function (pdf) with respect to Lebesgue measure of a distribution from a monotone polynomial tail (MPT) distribution family defined as follows.

\begin{definition}
A distribution with cdf $D$ and pdf $d$ is said to be an MPT distribution if $d$ is either an even function or a function whose support is the nonnegative real line, 
$d(x)x^{\gamma+1} \rightarrow C_d$ as $x\rightarrow \infty$ for some constant $\gamma >0$ and $C_d >0$, and for any $\theta > 1$, $d(x/\theta)/{d(x)}$ is a strictly increasing function in $x$ for $x > 0$.
\end{definition}

The study of decision procedures for distributions with general monotone likelihood ratio (MLR) properties dates back to \citet*{karlin1956theory}. The MLR property here ensures a simple form for the Bayes rule. These distributions have polynomial tails with $\gamma$, the polynomial tail heaviness index, specifying the heaviness of the tail. By L'Hospital's rule, one has $x^\gamma\{1-D(x)\} \rightarrow C_d/\gamma$. We focus on symmetric MPT distributions to make comparisons with normal distributions. For symmetric MPT, the scale family with pdf $d$ has MLR property in $|x|$. The MPT family includes many important distributions. Some examples are given in Table \ref{table:MPT_example}.

\begin{table}[htb]
\caption{Some common distributions in the MPT family.}{%
\begin{tabular}{lcc}
\hline
&  pdf & $C_d$ \\[1pt]
\hline
Student's t & $d(x) = C_d (x^2 + \gamma)^{-(\gamma+1)/2}, \ -\infty < x < \infty$ & $\frac{\gamma^{\gamma/2}\Gamma((\gamma+1)/2)}{{\pi}^{1/2}\Gamma(\gamma/2)}$ \\[5pt]
Pareto & $d(x) = C_d (x+1)^{-(\gamma+1)}, \ x > 0$ & $\gamma$\\[5pt]
Inverse Gamma & $d(x) = C_d x^{-\gamma-1}\exp(-1/x), \ x > 0$ & $1/\Gamma(\gamma)$\\[5pt]
\hline
\end{tabular}}
\label{table:MPT_example}
\end{table}

\begin{proposition}\label{prop:MPT_properties}
A distribution in the MPT family satisfies
\begin{enumerate}
\item \label{prop:1}
for $x>0$, $g(x) = d(x)x^{\gamma+1}$ is a strictly increasing function;
\item \label{prop:2}
for $x>0$, $h(x) = x^\gamma\{1-D(x)\}$ is a strictly increasing function;
\item \label{prop:3}
if $\theta > 1$, $\{1-D(x/\theta)\}/\{1-D(x)\}$ is a strictly increasing function for $x > 0$;
\item \label{prop:4}
for $x > 0$, $f(x) = 2xd(x)/\gamma + 2D(x)-1$ is a strictly increasing function.
\end{enumerate}
\end{proposition}

We assume that the common distribution of the $X$'s has cdf \begin{equation}\label{eq:mixture}
p D(x/\sigma_1) + (1-p)D(x/\sigma_0) \quad (0 <\sigma_0 < \sigma_1).
\end{equation}
The mixture has the distribution of $X_i$ determined by a latent Bernoulli random variable $s_i$ with success probability $p$, and the latent variables $s_1, \ldots, s_m$ are mutually independent. The cdf of $X_i$ is $D(x/\sigma_0)$ if $s_i=0$ and $D(x/\sigma_1)$ if $s_i=1$. We are interested in simultaneously testing $$H_{0i}: s_i=0~\textrm{vs.}~ H_{Ai}:s_i=1,~i=1, \ldots, m.$$
In a Bayesian hypothesis testing framework, the marginal distribution of the observations is often of the form (\ref{eq:mixture}). 
A specific example of our model is in stock selection, where $X_1, \ldots, X_m$ are the log returns of $m$ stocks, often modeled by Student's t distributions. We can equivalently assume $
X_i \, | \, \mu_i, \tau_i^2 \sim N(\mu_i, \tau_i^2)$ and $\tau_i^2 \sim \textsc{IG}(\gamma/2, \gamma/2)$. To test whether some stocks have extreme returns, if we assume $$\mu_i \, | \, \tau_i^2 \sim (1-p)N(0, \eta_{0}^2\tau_i^2) + pN(0, \eta_{1}^2\tau_i^2), ~ i=1,\ldots,m,$$ then the marginal distribution of $X_i$ is a mixture of two Student's t distributions with $\gamma$ degrees of freedom and $\sigma_j^2 = \eta_{j}^2 + 1, \, j=0,1$.

 Let $r_i : \mathbb{R}^m \rightarrow\{0,1\}$ be the decision rule used for the $i$th test. If ${r_i(\mathbf{X}) = 1}$, the null hypothesis is rejected; otherwise the null is not rejected. For each test, the loss is non-zero only if $r_i \neq s_i$, that is, only when a type \rom{1} or type \rom{2} error is made. Let $\delta_0$ and $\delta_A$ denote the respective losses of making a type \rom{1} and a type \rom{2} error. We assume that the overall loss of the $m$ tests is the sum of losses for individual tests. This additive loss structure is similar to the one in \citet*{lehmann1957theorya,lehmann1957theoryb} and in \citet*{bogdan2011asymptotic}.
To simplify matters, we take $$1+u=\frac{\sigma_1^2}{\sigma_0^2},\quad \delta=\frac{\delta_{0}}{\delta_{A}},\quad f=\frac{1-p}{p},\quad v=\delta f.$$ These parameters can vary with the number of tests $m$.

The Bayes risk of a multiple testing procedure is 
\begin{equation}\label{eq:risk}
R = p\delta_A\sum_{i=1}^m (vt_{1i} + t_{2i}),
\end{equation}
where $t_{1i}$ and $t_{2i}$ 
denote the probabilities of type \rom{1} and type \rom{2} errors for the $i$th test;
hereafter we call $vt_{1i}$ and $t_{2i}$ the type \rom{1} and type \rom{2} risk components, respectively.
The Bayes rule minimizing the Bayes risk can be shown to reject $H_{0i}$ if
$$\frac{{(1+u)^{-1/2}}d\left(({1+u})^{-1/2}{X_i}/{\sigma_0}\right)}{d(X_i/\sigma_0)} > v,~i=1,\ldots, m.$$
Due to the MLR property, the Bayes rule rejects $H_{0i}$ if
 $X^2_i/\sigma_0^2 \geq \omegaopt^2$, where $\omegaopt$ is the positive solution of the equation
 \begin{equation}\label{eq:oracle_threshold}
 \frac{(1+u)^{-1/2}d\left(\omegaopt(1+u)^{-1/2}\right)}{d(\omegaopt)}= v.
 \end{equation}
 As $p$ is unknown in practice, we call $\omegaopt$ the oracle threshold. 

We seek conditions for a multiple testing procedure to attain the Bayes oracle property under sparsity as $m \rightarrow \infty$. To impose sparsity, we assume $p \rightarrow 0$ as $m \rightarrow \infty$, and let $u \rightarrow \infty$ to ensure that the signals are strong enough to be discovered as $m \rightarrow \infty$. 
We assume that $\omegaopt^2/(1+u) \rightarrow C\in (0,\infty)$ to avoid having the power of an individual test going to zero or one. By (\ref{eq:oracle_threshold}) and the definition of MPT family, this is equivalent to assuming $vu^{-\gamma/2} \rightarrow C_0$, where 
\begin{equation}\label{eq:C0}
C_0=C^{(\gamma+1)/2}d(\sqrt{C})/C_d
\end{equation}
is a strictly increasing function in $C$, according to Proposition \ref{prop:MPT_properties}. If $\delta \rightarrow \delta_\infty=1$, this can be simplified to $pu^{\gamma/2} \rightarrow C_0^{-1}$. Intuitively, it is more difficult to distinguish between signals and noises if the data have fewer signals (smaller $p$) and smaller signal-to-noise ratio (smaller $u$). The assumption that $pu^{\gamma/2}$ converges to a constant guarantees the signals are identifiable, while the magnitude of the constant $C_0^{-1}$ indicates the intrinsic difficulty in identifying those signals. A larger $C$ reflects more difficulties and we call $C$ the difficulty index.

To summarize, we study the properties of multiple testing procedures under the asymptotic framework (as $m \rightarrow \infty$)
\begin{equation}\label{eq:framework}
p \rightarrow 0, \quad u \rightarrow \infty, \quad v{u^{-\gamma/{2}}} \rightarrow C_0,
\end{equation}
whereas the asymptotic framework in \citet*{bogdan2011asymptotic} has $p \rightarrow 0, u \rightarrow \infty, v \rightarrow \infty, 2\log (v)/u \rightarrow C$.
Noticing that the second and the third assumptions in \eqref{eq:framework} imply $v\rightarrow \infty$, the only difference between the two frameworks is the relation between $u$ and $v$. For normal distributions, the rate at which the signal strength $u$ increases to infinity is the logarithm of $v$,  while for polynomial-tailed distributions, it is a polynomial in $v$. 

\begin{proposition}\label{thm:oracle}
Let $C_1 = 2 \sqrt{C}d(\sqrt{C})/\gamma$ and ${C_2 = 2D(\sqrt{C})-1}$. Under (\ref{eq:framework}),
\begin{equation}\label{eq:oracle_asymptotic}
\omegaopt^2 \sim C (v / C_0)^{2/\gamma}, \quad t_{1i} = t_1 \sim C_1v^{-1}, \quad t_{2i} = t_2 \sim C_2,
\end{equation}
with $C_0$ as in (\ref{eq:C0}).
The corresponding Bayes risk is
\begin{equation}\label{eq:oracle_risk}
\Ropt \sim mp\delta_A(C_1 + C_2).
\end{equation}
\end{proposition}

By Proposition \ref{prop:MPT_properties}, $\Ropt$ is a strictly increasing function in $C$, which agrees with the interpretation of the difficulty index $C$; a more difficult multiple testing task leads to a higher Bayes risk.

Since $C_1$ and $C_2$ are the limits of type \rom{1} and type \rom{2} risk components of the oracle procedure, we call them asymptotically optimal type \rom{1} and type \rom{2} risk components. In \citet*{bogdan2011asymptotic} and \citet*{neuvial2012false}, the limiting Bayes risks of the oracle threshold are shown to depend solely on the type \rom{2} risk component. For polynomial-tailed distributions, neither risk component of the oracle procedure is negligible as the number of tests goes to infinity. In both models, the oracle probability of type \rom{1} errors goes to zero, but the probability decays at the rate of $v^{-1}$ for polynomial-tailed models, while for the normal model, the rate is faster. Besides the need for stronger signals to ensure detectability, this is yet another effect of heavy-tailed signals and noises.

\begin{definition}
A multiple testing rule is asymptotically Bayes optimal under sparsity (ABOS) under (\ref{eq:framework}) if its Bayes risk $R$ satisfies ${R}/{\Ropt} \rightarrow 1.$
\end{definition}

\par \bigskip

\setcounter{chapter}{3}
\setcounter{equation}{0} 
\noindent {\bf 3. Fixed and Random Thresholding Procedures}
\par \smallskip
In this section, we consider multiple testing procedures that reject the $i$th null hypothesis if $X_i^2/\sigma_0^2$ is greater than or equal to a threshold, which can be either non-data dependent (fixed) or data dependent (random). To distinguish them, we let $\omega^2$ denote the fixed threshold and $\omegahat^2$ denote the random threshold. 
For a fixed thresholding procedure, the events $\left\{r_i(\mathbf{X}) = 1 \mid s_i=0\right\}$ and $\left\{r_i(\mathbf{X}) = 0 \mid s_i=1 \right\}$ are based only on the $i$th observation $X_i$ with respective probabilities the same for each $i$. Therefore, the Bayes risk of a fixed thresholding procedure can be expressed as
$R = mp\delta_A(vt_1 + t_2),$ where $t_1 = 2\{1-D(\omega)\}$, and $t_2 = 2D(\omega(1+u)^{-1/2})-1.$ In contrast, the events $\left\{r_i(\mathbf{X}) = 1 \mid s_i=0\right\}$ and $\left\{r_i(\mathbf{X}) = 0 \mid s_i=1 \right\}$ for a random thresholding procedure are potentially based on all $m$ observations and the probabilities of type \rom{1} and type \rom{2} errors are not necessarily the same across different tests. 

\begin{theorem}\label{thm:fix_threshold}
A fixed thresholding multiple testing procedure that rejects $H_{0i}$ when ${{X^2_i}/{\sigma_0^2} \geq \omega^2, i=1, \ldots, m}$ is ABOS if and only if the threshold $\omega$ satisfies
\begin{equation}\label{eq:ABOS_threshold2}
\omega^2 / \omegaopt^2 \rightarrow 1,
\end{equation}
or, equivalently, with $C_0$ as in (\ref{eq:C0}),
\begin{equation}\label{eq:ABOS_threshold1}
\omega^2 = C (v/C_0)^{{2}/{\gamma}} (1+o(1)).
\end{equation}
\end{theorem}

It may appear that even if the type \rom{1} and type \rom{2} risk components do not tend to the corresponding asymptotically optimal risk components, there is still a chance that $R \sim \Ropt$. However, the proof shows that the two components have to converge to the corresponding optimal risk components individually in order to achieve ABOS. This observation is also true for the normal distribution, but as shown by Theorem 3.2 in \citet*{bogdan2011asymptotic}, two conditions, one for the type \rom{2} risk component, the other for the type \rom{1} risk component, are needed to guarantee this, while in our case, only one condition is required. In Remark 3.1 of \citet*{bogdan2011asymptotic}, the authors argued the reason for the extra condition for type \rom{1} error is that, for normal models, type \rom{1} errors are more sensitive to changes in the critical value than type \rom{2} errors. In their language, our Theorem \ref{thm:fix_threshold} shows that for polynomial-tailed models, type \rom{1} and type \rom{2} errors are equally sensitive to changes in the critical value.

\begin{theorem}\label{thm:random}
Under \eqref{eq:framework}, a random thresholding multiple testing procedure that rejects $H_{0i}$ if ${X_i^2}/{\sigma_0^2} \geq \omegahat^2$ is ABOS if for all $\epsilon > 0$,
\begin{equation}\label{eq:ABOS_random_type1}
\frac{1}{m}\sum_{i=1}^m P(|\omegahat-\omegaopt| > \epsilon v^{1/\gamma} | s_i = 0) = o(v^{-1}),
\end{equation}
\begin{equation}\label{eq:ABOS_random_type2}
\frac{1}{m}\sum_{i=1}^m P(|\omegahat-\omegaopt| > \epsilon v^{1/\gamma} | s_i = 1) = o(1).
\end{equation}
If $\delta$ does not converge to zero as $m \rightarrow \infty$, then a random thresholding procedure is ABOS if for all $\epsilon > 0$,
\begin{equation}\label{eq:ABOS_random2}
P(|\omegahat/\omegaopt - 1| > \epsilon) = o(v^{-1}).
\end{equation}
\end{theorem}

An equivalent condition to \eqref{eq:ABOS_random2} is that
$P(|\omegahat-\omegaopt| > \epsilon v^{1/\gamma}) = o(v^{-1}).$
Theorem \ref{thm:random} continues to hold if the oracle threshold is replaced by the threshold of a fixed thresholding procedure that is ABOS. The left hand sides of \eqref{eq:ABOS_random_type1} and \eqref{eq:ABOS_random_type2} can be interpreted as the average departures of the probabilities of type \rom{1} and type \rom{2} errors of a random thresholding procedure from the corresponding errors of the Bayes oracle. 
Although less general than \eqref{eq:ABOS_random_type1} and \eqref{eq:ABOS_random_type2}, condition \eqref{eq:ABOS_random2} could be easier to verify in practice because of the symmetry of the distribution of $m$ observations. 


As implied by Theorem \ref{thm:fix_threshold}, to obtain an ABOS fixed thresholding procedure, the fixed threshold itself is very likely to contain unknown parameters. In contrast, a random threshold consists of observed data only. For example, it could appear as an estimator of an ABOS fixed threshold. Therefore, a random thresholding procedure is naturally an implementable procedure, and, in this sense, Theorem \ref{thm:random} provides a more practical result. 
\par \bigskip

\setcounter{chapter}{4}
\setcounter{equation}{0} 
\noindent {\bf 4. ABOS of Several Special Procedures}
\par \smallskip
\noindent {\bf 4.1 Procedures controlling BFDR}
\par \smallskip
\citet*{benjamini1995controlling} introduced \textsc{FDR} as a less stringent error measure than the familywise error rate, $\textsc{FDR} = E\left({V}/{R}\right),$ where $R$ is the number of total rejections and $V$ is the number of false rejections. \citet*{storey2003positive} argued that the positive false discovery rate (\textsc{pFDR}), $\textsc{pFDR} = E\left({V}/{R} \mid R >0 \right),$ can overcome some of the concerns in \citet*{benjamini1995controlling} and, under certain conditions, it coincides with the Bayesian false discovery rate (\textsc{BFDR}) of \citet*{efron2002empirical},
\begin{equation*}
\textsc{BFDR} = P(H_{0i}~\textrm{is true} \mid H_{0i}~\textrm{was rejected}) = \frac{(1-p)t_1}{(1-p)t_1+p(1-t_2)}.
\end{equation*}
For a fixed thresholding procedure, the threshold $\omega$ and the \textsc{BFDR} level $\alpha$ are linked by
\begin{equation}\label{eq:alpha_BFDR}
\frac{(1-p)\{1-D(\omega)\}}{(1-p)\{1-D(\omega)\} + p\left\{1-D({\omega}{{(1+u)}^{-1/2}})\right\}} = \alpha,
\end{equation}
or equivalently,
\begin{equation}\label{eq:def_BFDR2}
\frac{1-D(\omega)}{1-D\left({\omega}{{(1+u)^{-1/2}}}\right)} = \frac{r_\alpha}{f},
\end{equation}
where $r_\alpha={\alpha}/{(1-\alpha)}$.
Since we have already found a necessary and sufficient condition for the fixed thresholding procedure to be ABOS, by using (\ref{eq:def_BFDR2}), we are able to find conditions on $\alpha$ (depending on $m$) such that the BFDR controlling procedure is ABOS.

An alternative expression of the \textsc{BFDR} level $\alpha$ in \eqref{eq:alpha_BFDR} is that $\alpha = (1-p)\{1-p+pg(\omega)\}^{-1},$ where $g(\omega) = \{1-D(\omega)\}^{-1}\left\{1-D({\omega}{{(1+u)}^{-1/2}})\right\}$. Property \ref{prop:3} of Proposition \ref{prop:MPT_properties} shows that $g$ is a strictly increasing function in $\omega$, and,  since $1-D(x) \sim x^{-\gamma}$ as $x\rightarrow \infty$, $g(\omega) \rightarrow (1+u)^{\gamma/2}$ as $\omega \rightarrow \infty$. As $g(\omega)=1$ if $\omega=0$, the \textsc{BFDR} of a finite fixed threshold procedure for a given $m$ can only be controlled within the interval $I = (\beta^*, 1-p \, ]$, where 
\begin{equation}\label{eq:beta_our}
\beta^*=\{1 + (1+u)^{\gamma/2}/f\}^{-1}.
\end{equation}
Thus with a \textsc{BFDR} level less than $\beta^*$, the fixed threshold has to be infinite. In this case, none of the $m$ null hypotheses is rejected and the power of an individual test is zero. Since our asymptotic framework requires the power of an individual test to go to a nonzero constant, we confine $\alpha$ to the interval $I$.

To distinguish from general fixed thresholds, $\omega_{B_\alpha}^2$ is used to denote the fixed threshold controlling \textsc{BFDR} under $\alpha$, and the subscript $\alpha$ is omitted if there is no ambiguity.

\begin{proposition}\label{thm:BFDR}
A fixed thresholding procedure controlling \textsc{BFDR} under $\alpha$ is ABOS if and only if
\begin{equation}\label{eq:BFDR_condition}
{\delta r_\alpha} \rightarrow {C_1}/{(1-C_2)}.
\end{equation}
The threshold is of the form
\begin{equation}\label{eq:BFDR_threshold}
\omega^2_B = C_B\left(\frac{f}{r_\alpha}\right)^{{2}/{\gamma}}(1+o(1)),
\end{equation}
where $C_B = \left[(C_d/\gamma)/\{1-D(\sqrt{C})\}\right]^{{2}/{\gamma}}$ and $C_1, C_2$ as in Proposition \ref{thm:oracle}.
\end{proposition}

Condition (\ref{eq:BFDR_condition}) implies that if either one of $\delta$ and $\alpha$ goes to a positive constant, the other is forced to converge to a positive constant as well. For example, if $\delta$ converges to a positive constant $\delta_{\infty}$, then $\alpha \rightarrow \alpha_{\infty}$ where $\alpha_\infty$ is defined by 
\begin{equation}\label{eq:alpha_infty}
\alpha_\infty = \frac{1}{1+\delta_\infty(1-C_2)/C_1}.
\end{equation}
Also, as more penalty is imposed for type \rom{2} errors ($\delta \rightarrow 0$) as $m \rightarrow \infty$, then no control is made on \textsc{BFDR} since $\alpha \rightarrow 1$ as $m \rightarrow \infty$. 
\par \bigskip

\noindent {\bf 4.2 Genovese-Wasserman and Benjamini-Hochberg Procedures}
\par \smallskip
Let $Z_i = \left| {X_i}/{\sigma_0}\right|$ and $p_i = 2\{1-D(Z_i)\}$ denote the p-values for the $m$ tests, ordered as $p_{(1)} \leq \cdots \leq p_{(m)}$. The Benjamini-Hochberg procedure at FDR level $\alpha$ then looks for the largest $k$, denoted by $\hat k$, that satisfies $p_{(k)}\leq \alpha k/m$ and rejects all the tests whose p-values are less than or equal to $p_{(\hat k)}$. This is equivalent to rejecting the null hypothesis $H_{0i}$ if $Z_i^2 \geq \omegaBH^2$, where
\begin{equation}\label{eq:BH_threshold}
\omegaBH = \inf\left\{ y : \frac{2\{1-D(y)\}}{1-\hat{F}(y)} \leq \alpha \right\},
\end{equation}
$F$ being the common cdf of $Z_i$'s and $1-\hat F(y) = \# \{Z_i \geq y\}/m$. Thus the Benjamini-Hochberg procedure is  a random thresholding procedure. To study the ABOS of the Benjamini-Hochberg procedure via Theorem \ref{thm:random}, we need to compare (\ref{eq:BH_threshold}) with a fixed ABOS threshold. \citet*{genovese2002operating} showed that the Benjamini-Hochberg procedure can be approximated by a fixed thresholding procedure whose threshold $\omegaGW$ is the solution of
\begin{equation}\label{eq:GW_def}
\frac{1-D(\omegaGW)}{(1-p)\{1-D(\omegaGW)\}+p\{1-D(\omegaGW{(1+u)^{-1/2}})\}}=\alpha. 
\end{equation}

\begin{proposition}\label{thm:GW}
If $\alpha \not\rightarrow 1$, the rule that rejects the null hypothesis $H_{0i}$ when ${X_i^2}/{\sigma_0^2} \geq \omegaGW^2$ is ABOS if and only if (\ref{eq:BFDR_condition}) holds. In this case, with $C_B$ as in Proposition \ref{thm:BFDR},
$$\omegaGW^2 = C_B\left(\frac{f}{r_\alpha}\right)^{{2}/{\gamma}}(1+o(1)).$$
\end{proposition}

\begin{theorem}\label{thm:BH}
If
\begin{equation}\label{eq:mp1}
p \propto 1/\log(m) ~ \textrm{or}~p \propto m^{-\kappa}~\textrm{for some}~0 <\kappa < 1,
\end{equation}
\begin{equation}\label{eq:delta}
\delta \rightarrow {\delta_\infty} > 0,
\end{equation}
the Benjamini-Hochberg procedure at \textsc{FDR} level $\alpha$ is ABOS if $\alpha \rightarrow \alpha_\infty$, where $\alpha_\infty$ is at \eqref{eq:alpha_infty}.
\end{theorem}

The oracle Bayes rule balances type \rom{1} and type \rom{2} errors with the consideration of loss for each type of error. 
The optimal \textsc{FDR} level given in (\ref{eq:alpha_infty}) is indeed the result of balancing since it is determined by the limiting loss ratio $\delta_\infty$ and the asymptotically optimal risk components $C_1, C_2$. 


The asymptotically optimal \textsc{FDR} level depends on the difficulty index $C$, which is usually an unknown parameter. Although not having a conclusive answer, we discuss how to find a practically usable \textsc{FDR} level in Section 6.
\par \bigskip

\setcounter{chapter}{5}
\setcounter{equation}{0} 
\noindent {\bf 5. Simulation Results}\label{sec:simulation}
\par \smallskip
We compared the performances of the Bayes oracle, the Benjamini-Hochberg procedure with the optimal \textsc{FDR} level, and the Benjamini-Hochberg procedure with \textsc{FDR} level $1/\log(m)$, through simulation studies. The \textsc{FDR} level $1/\log(m)$ was proved to be ABOS for the normal distributions by \citet*{bogdan2011asymptotic}. The simulation study in \citet*{ghosh2013asymptotic} demonstrated its effectiveness in producing a misclassification probability curve similar to the one obtained from the oracle procedure. We considered this \textsc{FDR} level to illustrate the consequence of applying a multiple testing procedure regardless of the underlying distribution. We write the $\alpha$-BH procedure for the Benjamini-Hochberg procedure with \textsc{FDR} level $\alpha$.

The comparison of performances was done under two scenarios. In the first, we took the sparsity parameter to vary with number of tests $m$, $p=m^{-0.5}$. We recorded the risks of the Bayes oracle and Benjamini-Hochberg procedures with different \textsc{FDR} levels. In the second scenario, with $m=10^6$, we examined the behavior of the multiple testing procedures with changing values of $p$. In both scenarios, we fixed the parameters of the loss function, $\delta_A$ and $\delta_0$, to be 1. We considered combinations of polynomial tail heaviness index $\gamma$ and the difficulty index $C$, choosing from $\{3,10\}$ and $\{0.1, 1, 10\}$, respectively. For each combination, 1000 data sets of $X_i, i=1, \ldots, m$ were generated from the mixture distribution (\ref{eq:mixture}) with $\sigma_0^2=1$, $\sigma_1^2 = \sigma_0^2 (1+u)$, where $u=(v/C_0)^{2/\gamma}$. Pareto and Student's t distributions in the MPT distribution family were considered in the simulation.
\par \bigskip

\noindent {\bf 5.1 Results from scenario 1}
\par
\begin{figure}[htb]
\centering
	\begin{subfigure}[b]{0.32\textwidth}
		\caption{$m=10^2$ Risk Ratios}
		\includegraphics[width=\textwidth]{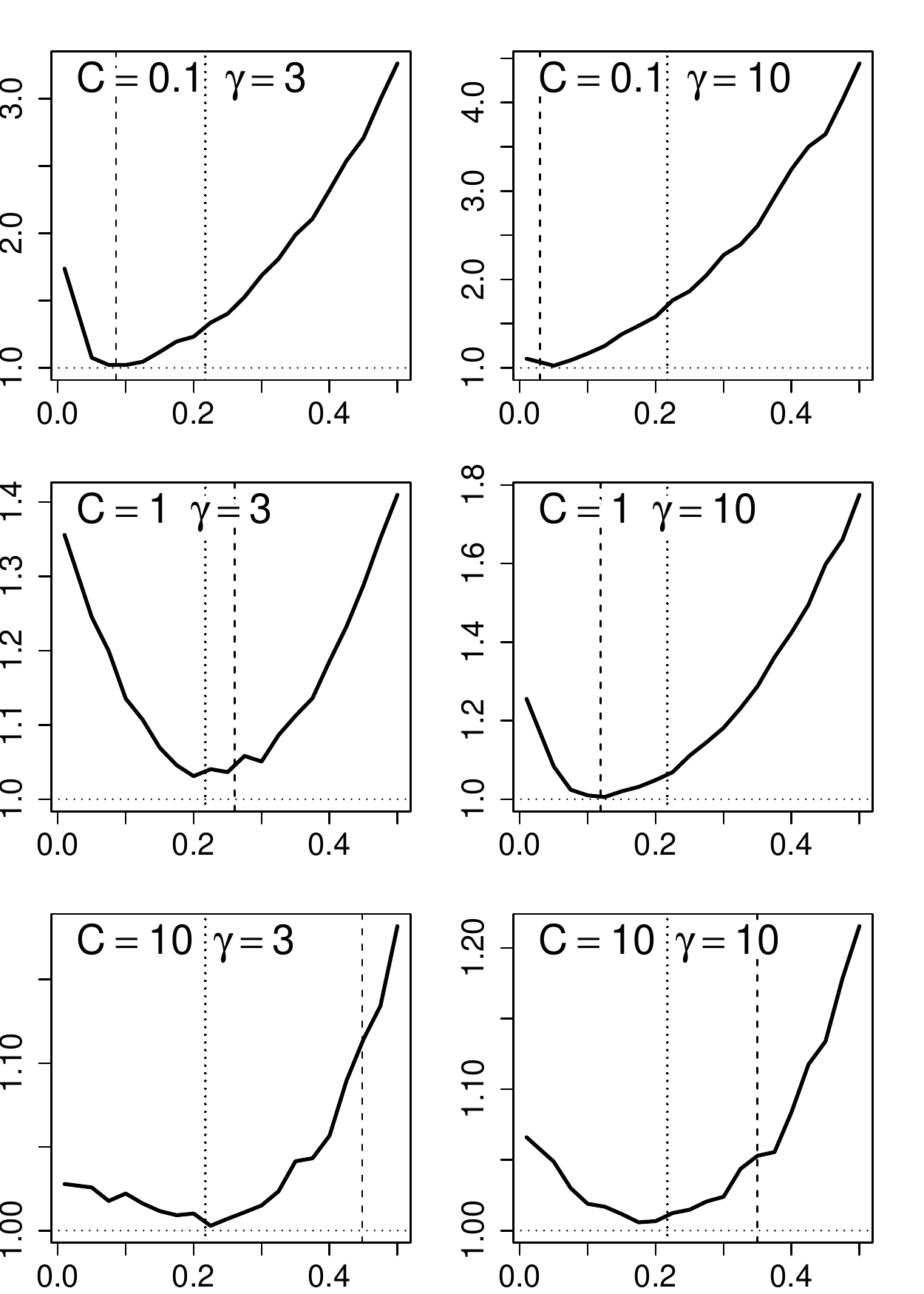}
	\end{subfigure}
	\,
	\begin{subfigure}[b]{0.32\textwidth}
		\caption{$m=10^6$ Risk Ratios}
		\includegraphics[width=\textwidth]{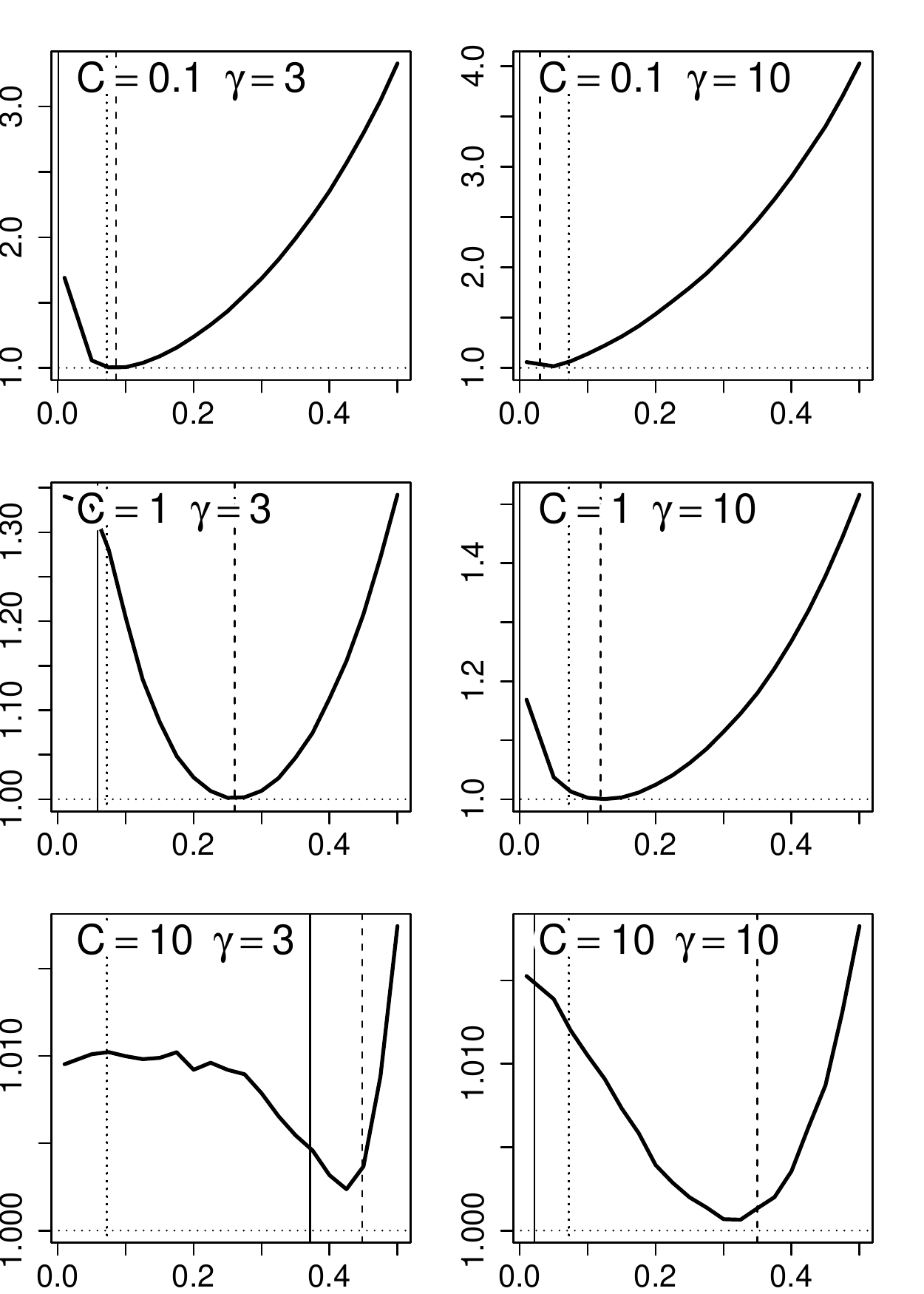}	
	\end{subfigure}
	\,
	\begin{subfigure}[b]{0.32\textwidth}
		\caption{$m=10^6$ P2}
		\includegraphics[width=\textwidth]{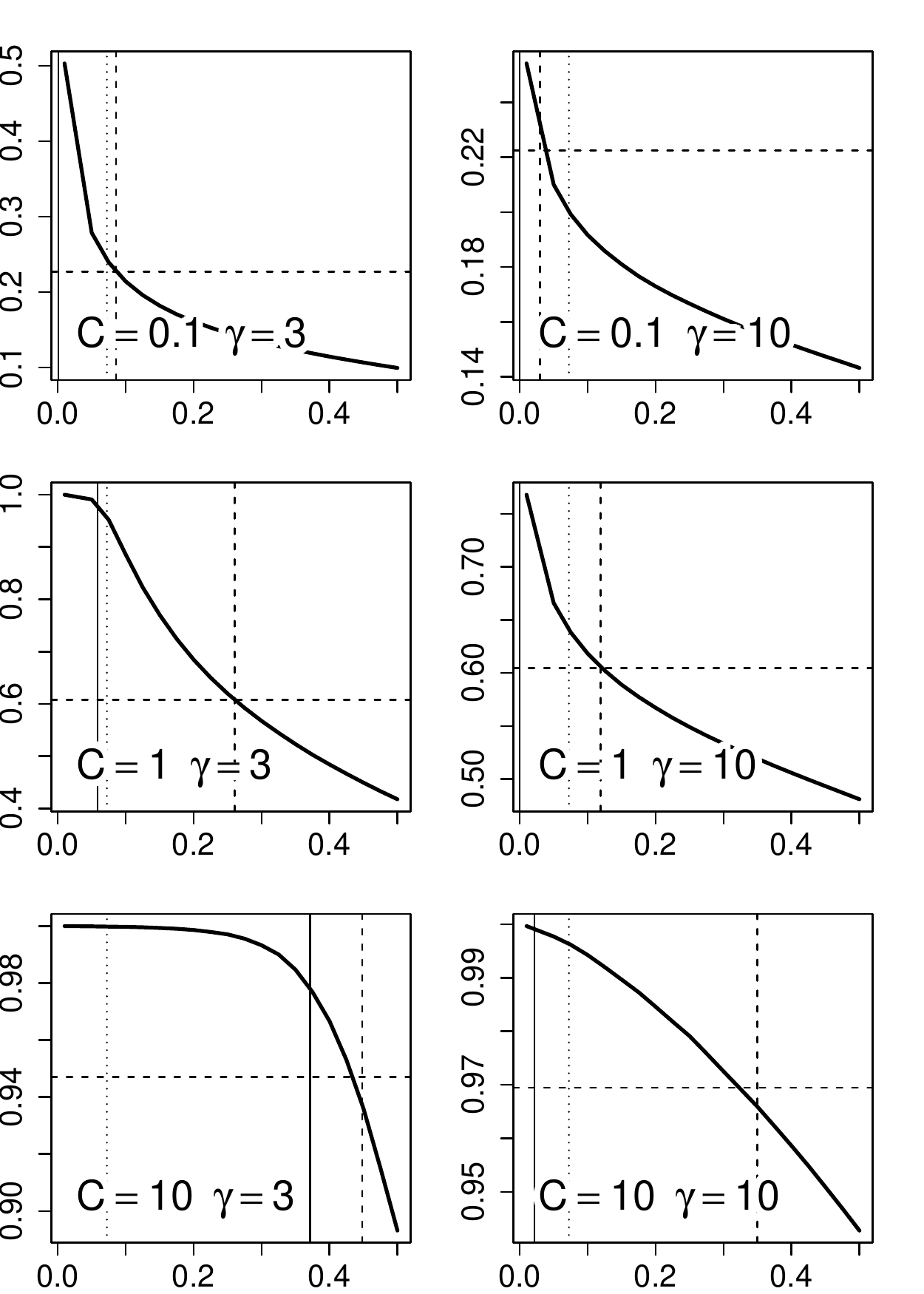}	
	\end{subfigure}
\caption{Plots of Bayes risk ratios (panels (a) and (b)) and probabilities of type \rom{2} error (P2; panel (c)) against \textsc{FDR} level $\alpha$ for Student's T distributions. The dashed vertical lines denote the asymptotically optimal $\alpha$ and the dotted vertical lines denote $1/\log(m)$. The solid vertical lines in panels (b) and (c) denote $\beta^*_\infty$ defined in \eqref{eq:betastar_inf}. The dashed horizontal lines in panel (c) represent P2 of the Bayes oracle. The ranges of vertical axes are different across plots.}\label{fig:t}
\end{figure}

\begin{figure}[htb]
\centering
	\begin{subfigure}[b]{0.32\textwidth}
		\caption{$m=10^2$}
		\includegraphics[width=\textwidth]{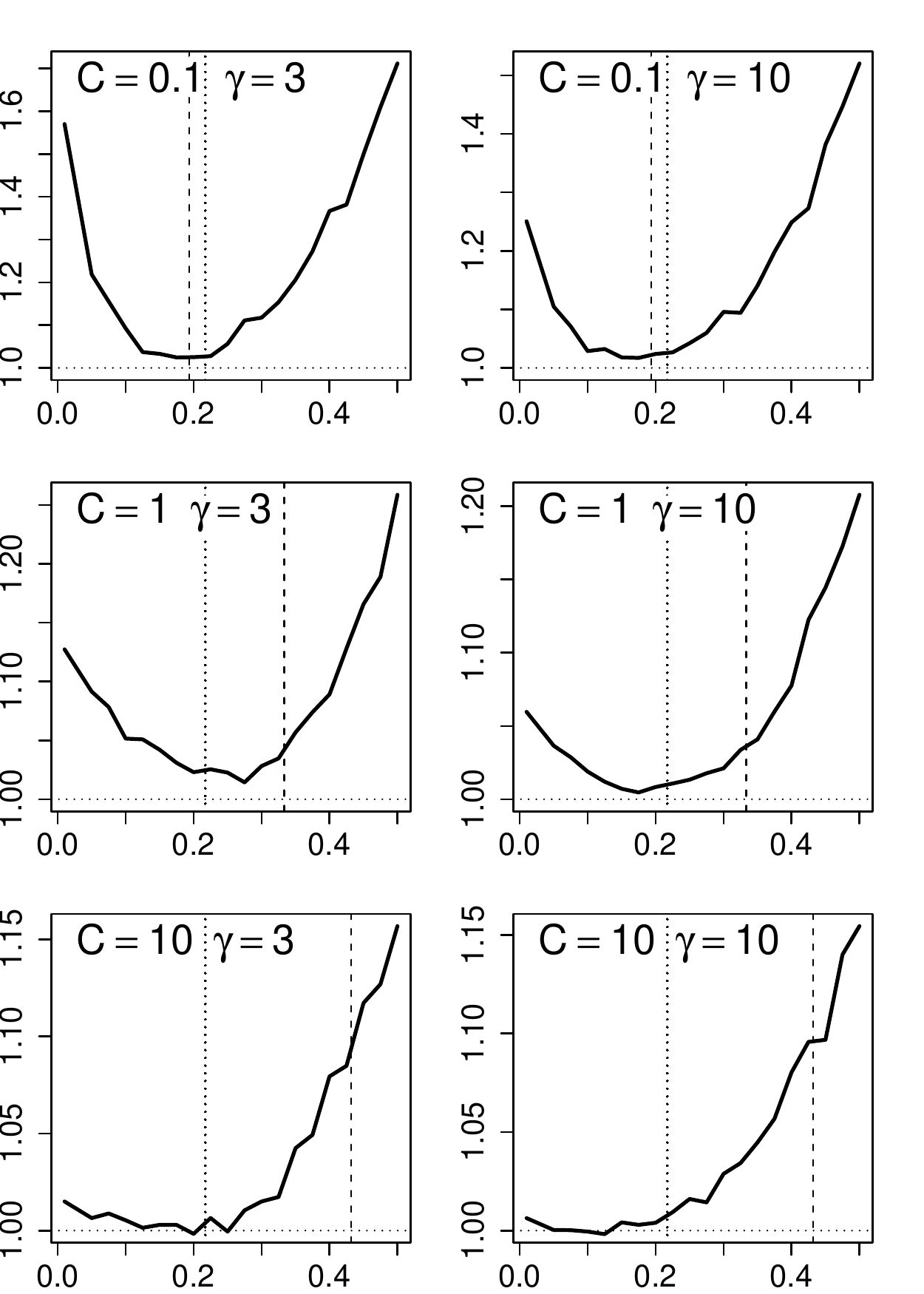}
	\end{subfigure}
	\,
	\begin{subfigure}[b]{0.32\textwidth}
		\caption{$m=10^6$}
		\includegraphics[width=\textwidth]{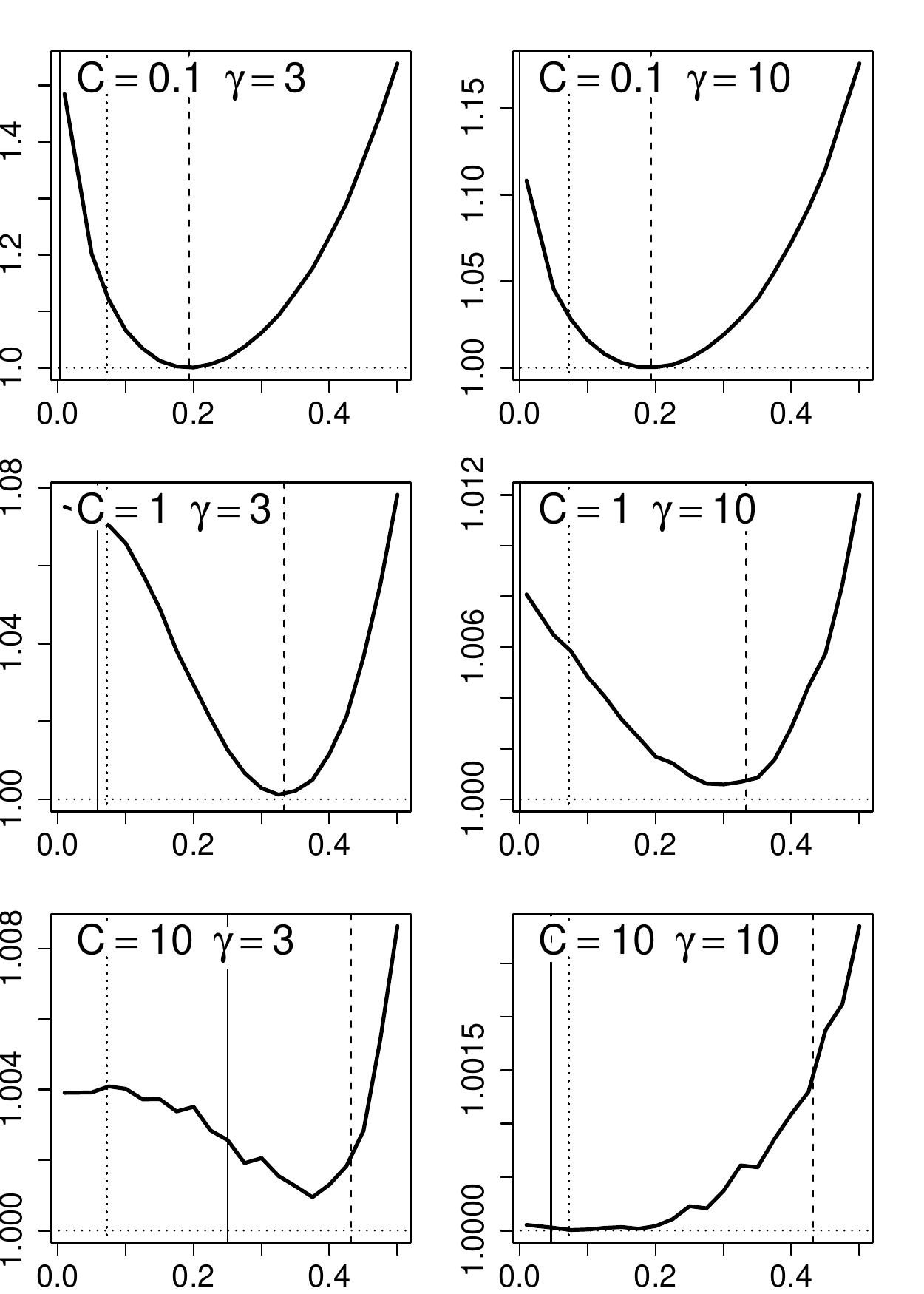}	
	\end{subfigure}
	\,
	\begin{subfigure}[b]{0.32\textwidth}
		\caption{$m=10^6$ P2}
		\includegraphics[width=\textwidth]{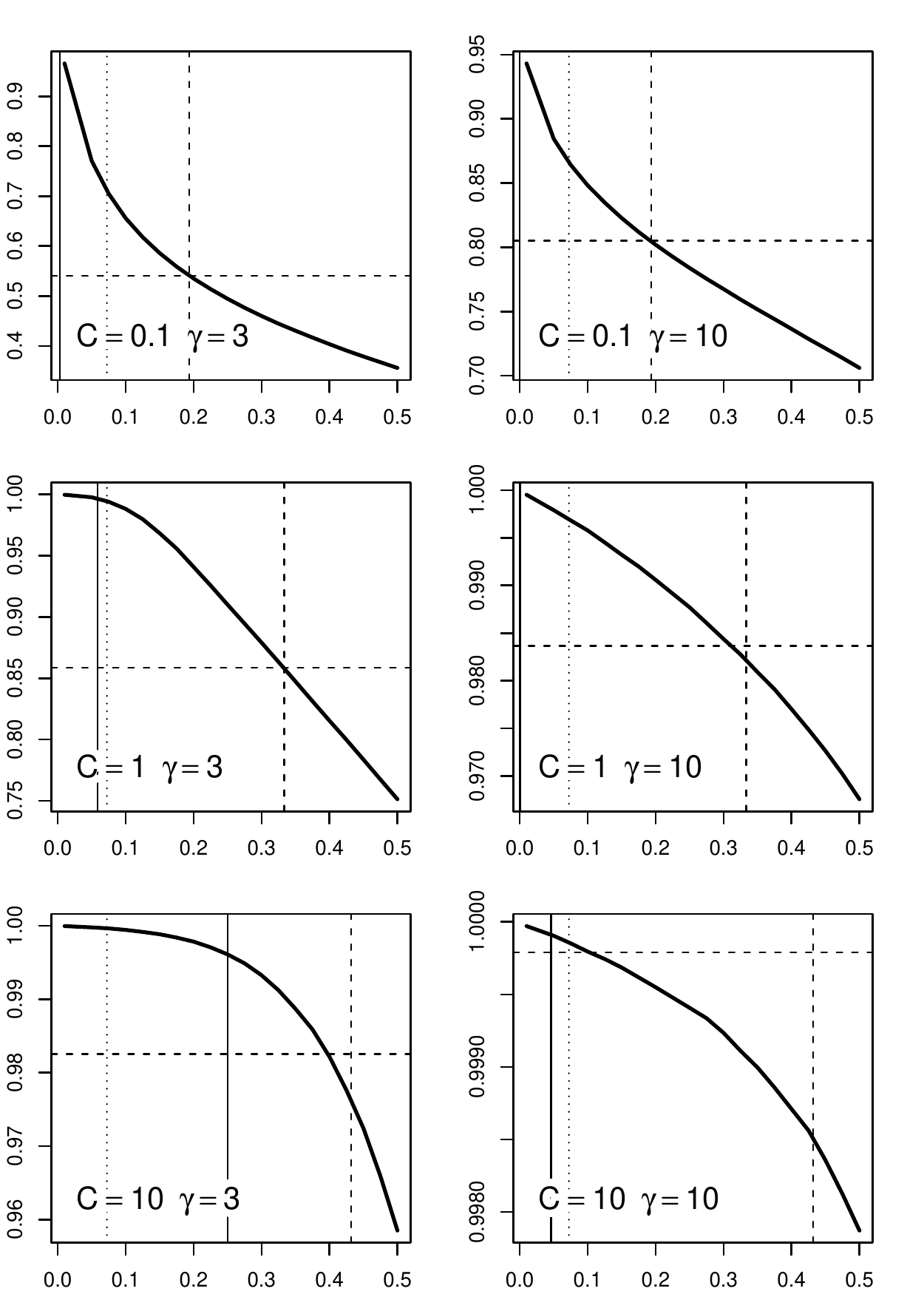}	
	\end{subfigure}
\caption{Plots of Bayes risk ratios (panels (a) and (b)) and probabilities of type \rom{2} error (P2; panel (c)) against \textsc{FDR} level $\alpha$ for Pareto distributions. The dashed vertical lines denote the asymptotically optimal $\alpha$ and the dotted vertical lines denote $1/\log(m)$. The solid vertical lines in panels (b) and (c) denote $\beta^*_\infty$ defined in \eqref{eq:betastar_inf}. The dashed horizontal lines in panel (c) represent P2 of the Bayes oracle. The ranges of vertical axes are different across plots.}\label{fig:pareto}
\end{figure}

 \smallskip
The average risks based on 1000 replicates were used to estimate the Bayes risks of the two procedures and to find the Bayes risk ratio of the BH procedure to the oracle. Panels (a) and (b) of Figures \ref{fig:t} and \ref{fig:pareto} show the plots of Bayes risk ratios against \textsc{FDR} level $\alpha$. In the plots, the dashed vertical lines denote the asymptotically optimal \textsc{FDR} level $\alpha_\infty$ as at (\ref{eq:alpha_infty}). When $m=10^6$, the risk ratios at $\alpha=\alpha_\infty$ are close to 1 and almost reach the lowest point of the curve. When $m=10^2$, the risk ratios at $\alpha=\alpha_\infty$ is not as close to the minimum as in the case $m=10^6$, but the deviations are still moderate. This observation does not conflict with the asymptotic results we have established, but it suggests that the study of non-asymptotic results or the convergence rate of asymptotic results may be helpful to find a better $\alpha$ for smaller $m$. In the figure, the dotted vertical lines in the plots are $\alpha = 1/\log(m)$. In some situations, this choice of \textsc{FDR} level does a better job than $\alpha_\infty$, but it can also lead to a risk ratio away from 1 in other situations.

In the plots, the range of the Bayes risk ratios is narrower for a larger $C$. This is probably because the denominator of the ratio, the oracle Bayes risk, is an increasing function in $C$. 

The optimal \textsc{FDR} level $\alpha_\infty$ (dashed vertical lines in Figures \ref{fig:t} and \ref{fig:pareto}) increases as the difficulty index $C$ increases. With a larger $C$, which signifies more difficulties in identifying signals from noises, the \textsc{FDR} can only be controlled at a higher level to achieve asymptotic Bayesian optimality. For both Student's t and Pareto distributions, when $C=10$, $\alpha_\infty$ is close to 0.5, which could hardly provide satisfactory control of false discoveries in practice.

\par \bigskip

\noindent {\bf 5.2 Results from scenario 2}
\par
\begin{figure}[htb]
\centering
	\begin{subfigure}[b]{0.32\textwidth}
		\caption{MP}
		\includegraphics[width=\textwidth]{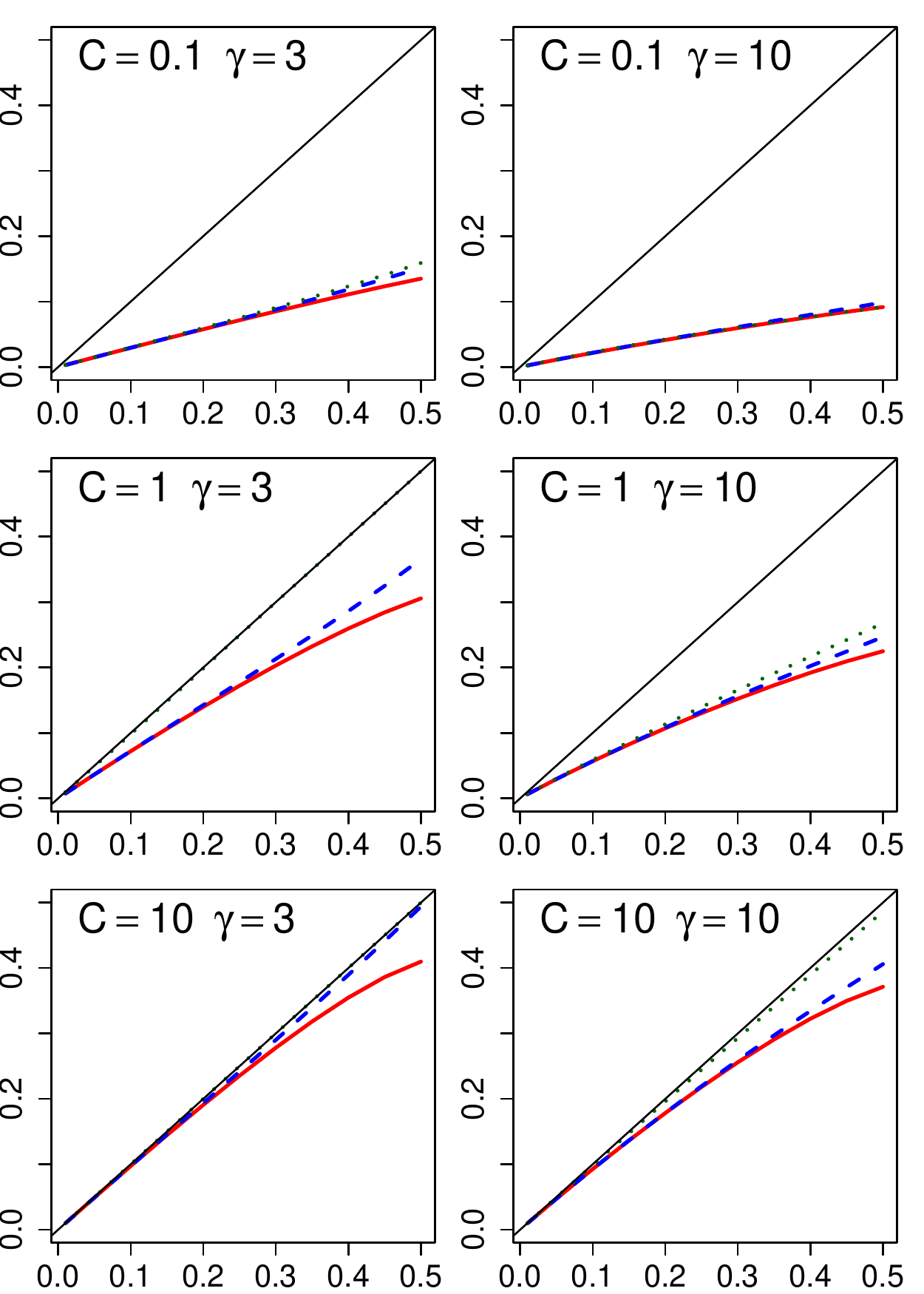}
	\end{subfigure}
	\,
	\begin{subfigure}[b]{0.32\textwidth}
		\caption{P1}
		\includegraphics[width=\textwidth]{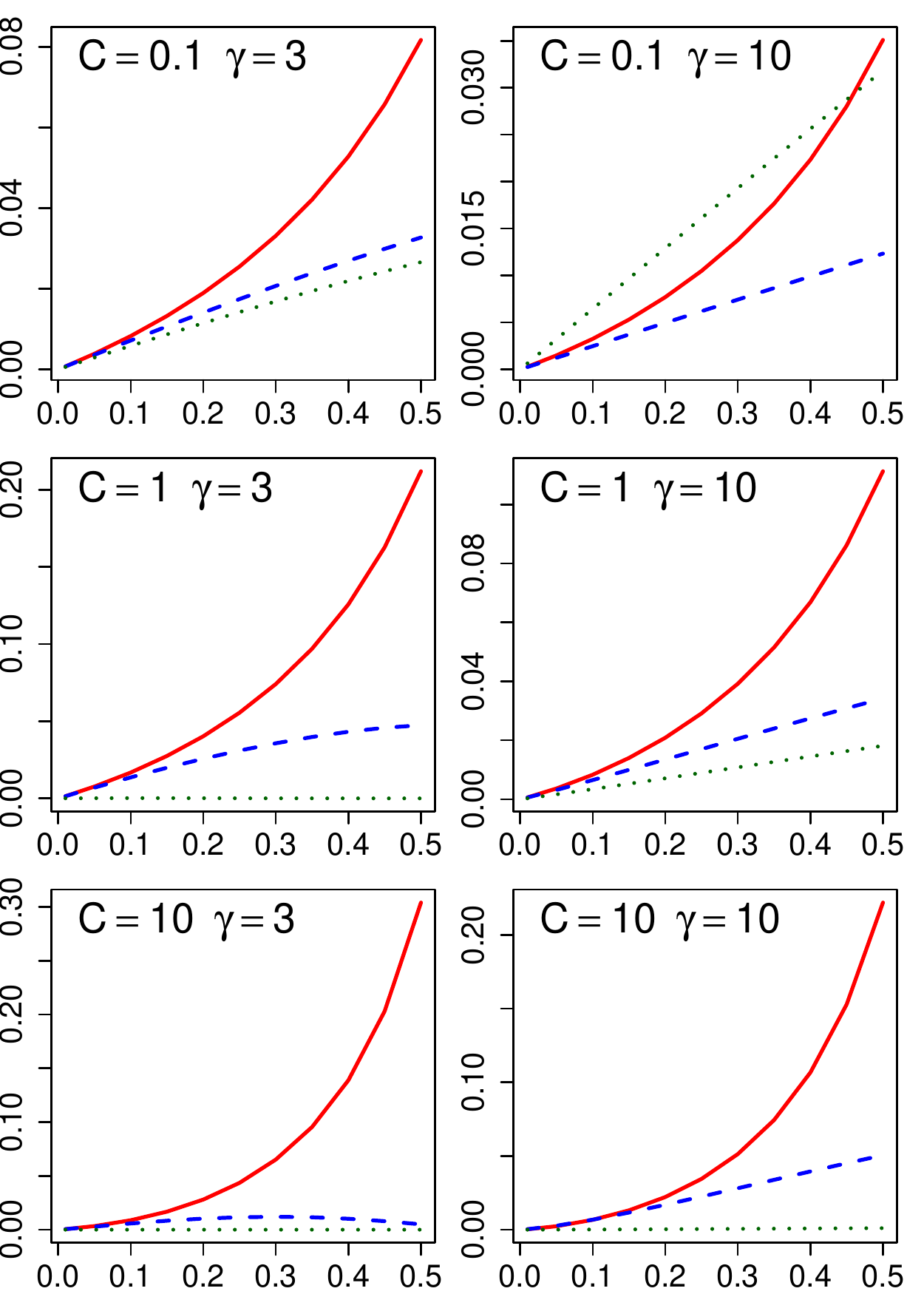}	
	\end{subfigure}
	\,
	\begin{subfigure}[b]{0.32\textwidth}
		\caption{P2}
		\includegraphics[width=\textwidth]{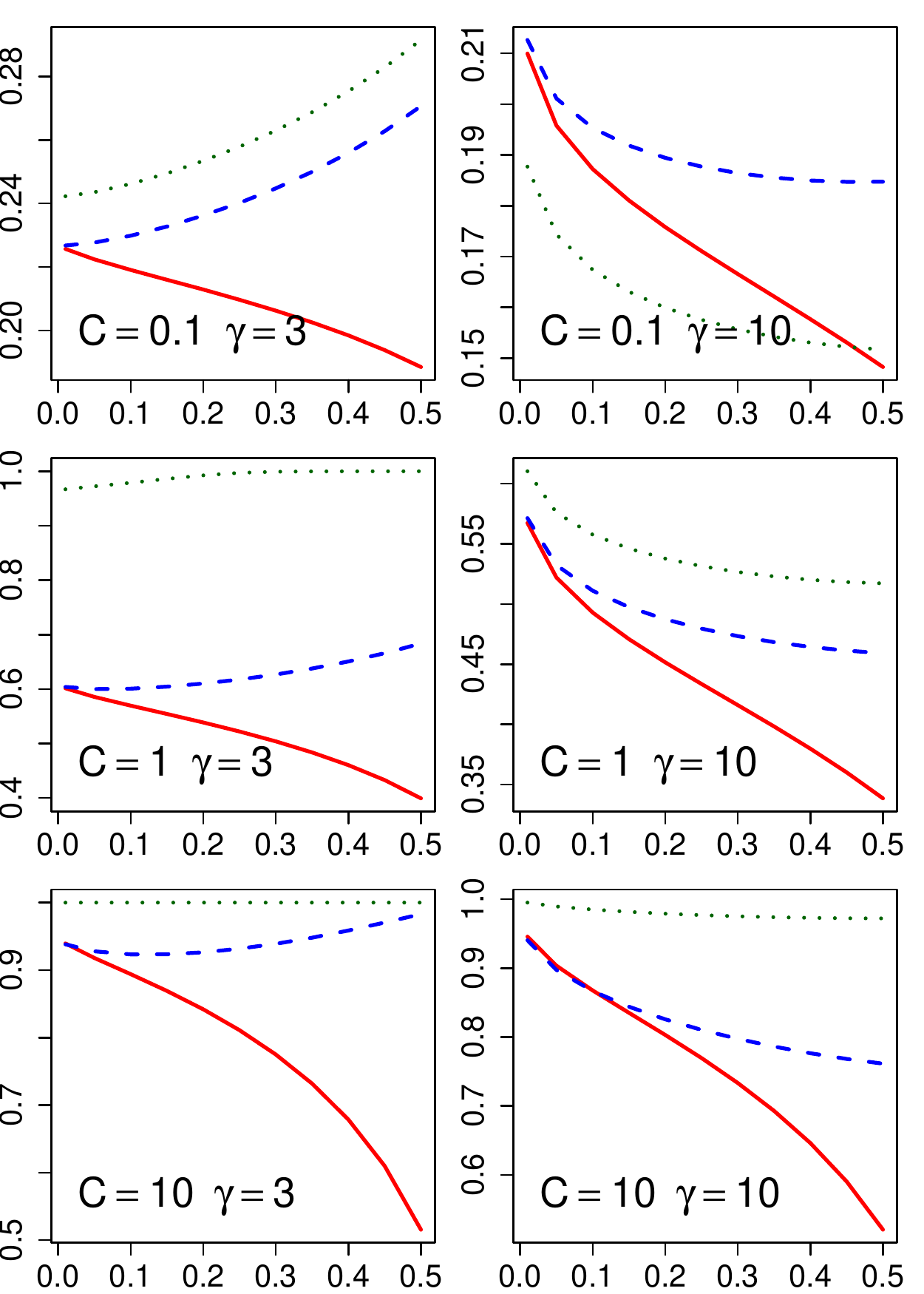}	
	\end{subfigure}
\caption{Plots of MP, P1, and P2 against sparsity parameter $p$ for Student's T distributions. The red solid, blue dashed, and green dotted lines represent the Bayes oracle, the $\alpha_\infty$-BH procedure and the $1/\log(m)$-BH procedure respectively. In panel (a), the diagonal lines are the MP=$p$ line. The ranges of the vertical axes are different across the plots of P1 and P2.}\label{fig:t2}
\end{figure}

\begin{figure}[htb]
\centering
	\begin{subfigure}[b]{0.32\textwidth}
		\caption{MP}
		\includegraphics[width=\textwidth]{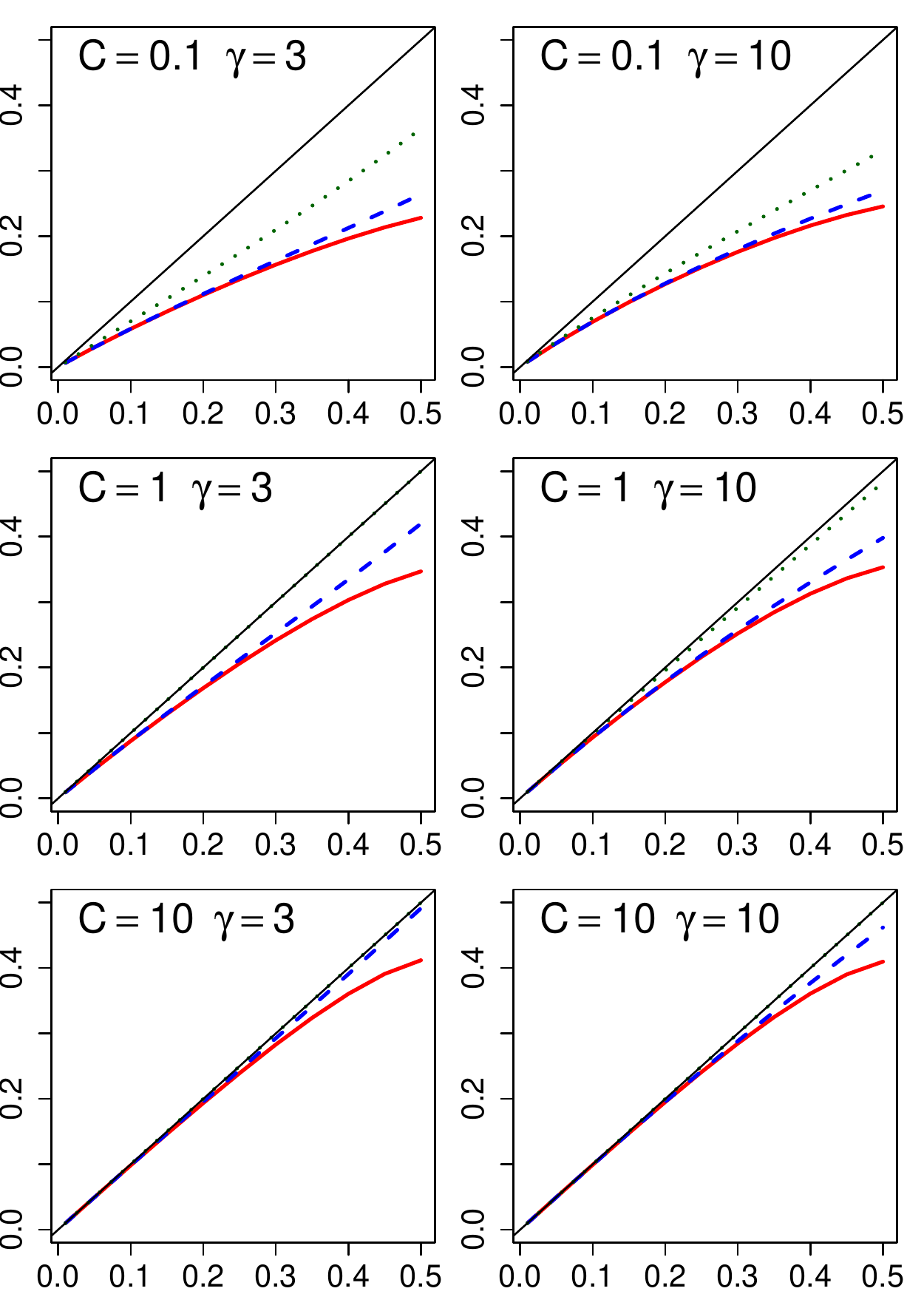}
	\end{subfigure}
	\,
	\begin{subfigure}[b]{0.32\textwidth}
		\caption{P1}
		\includegraphics[width=\textwidth]{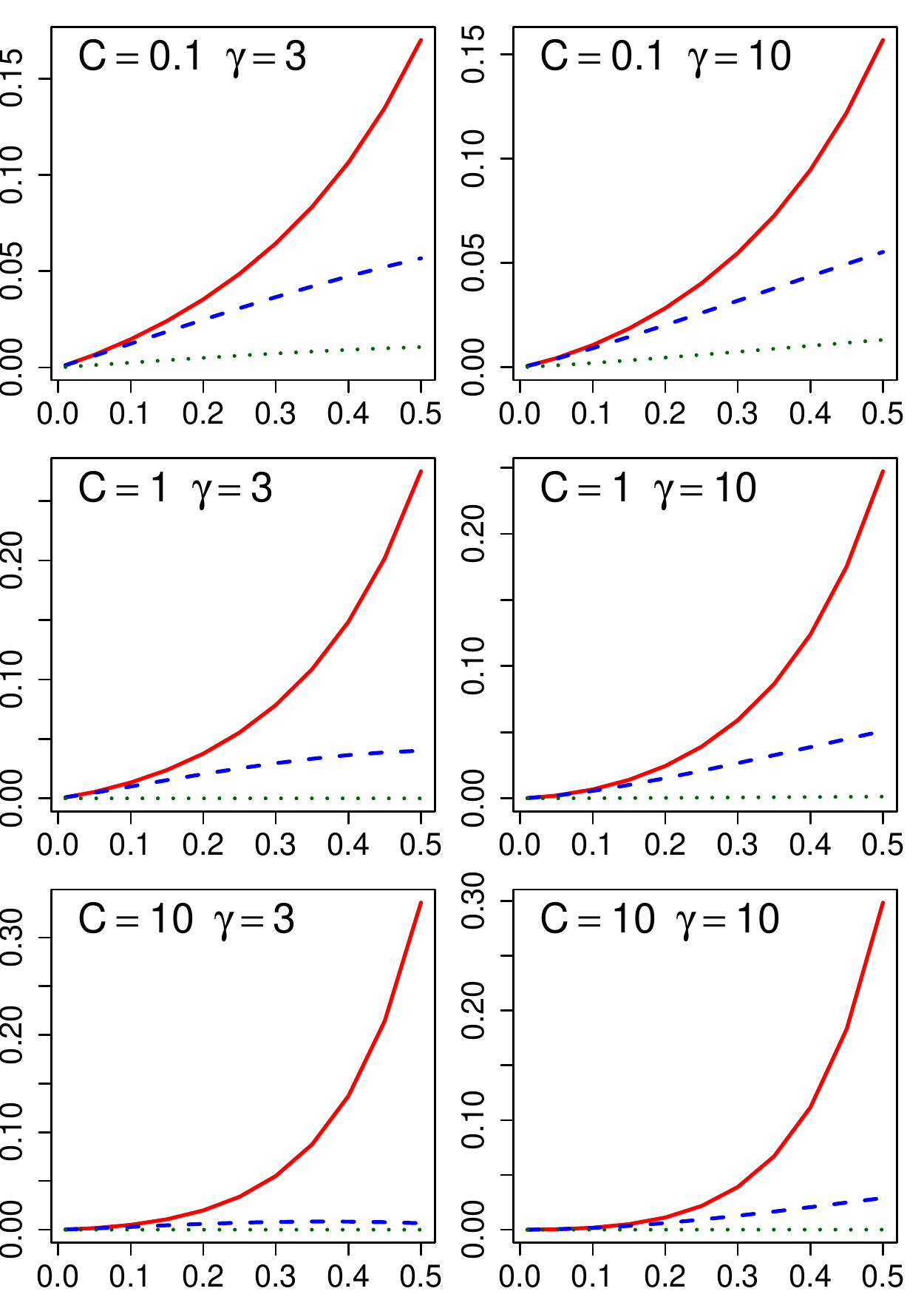}	
	\end{subfigure}
	\,
	\begin{subfigure}[b]{0.32\textwidth}
		\caption{P2}
		\includegraphics[width=\textwidth]{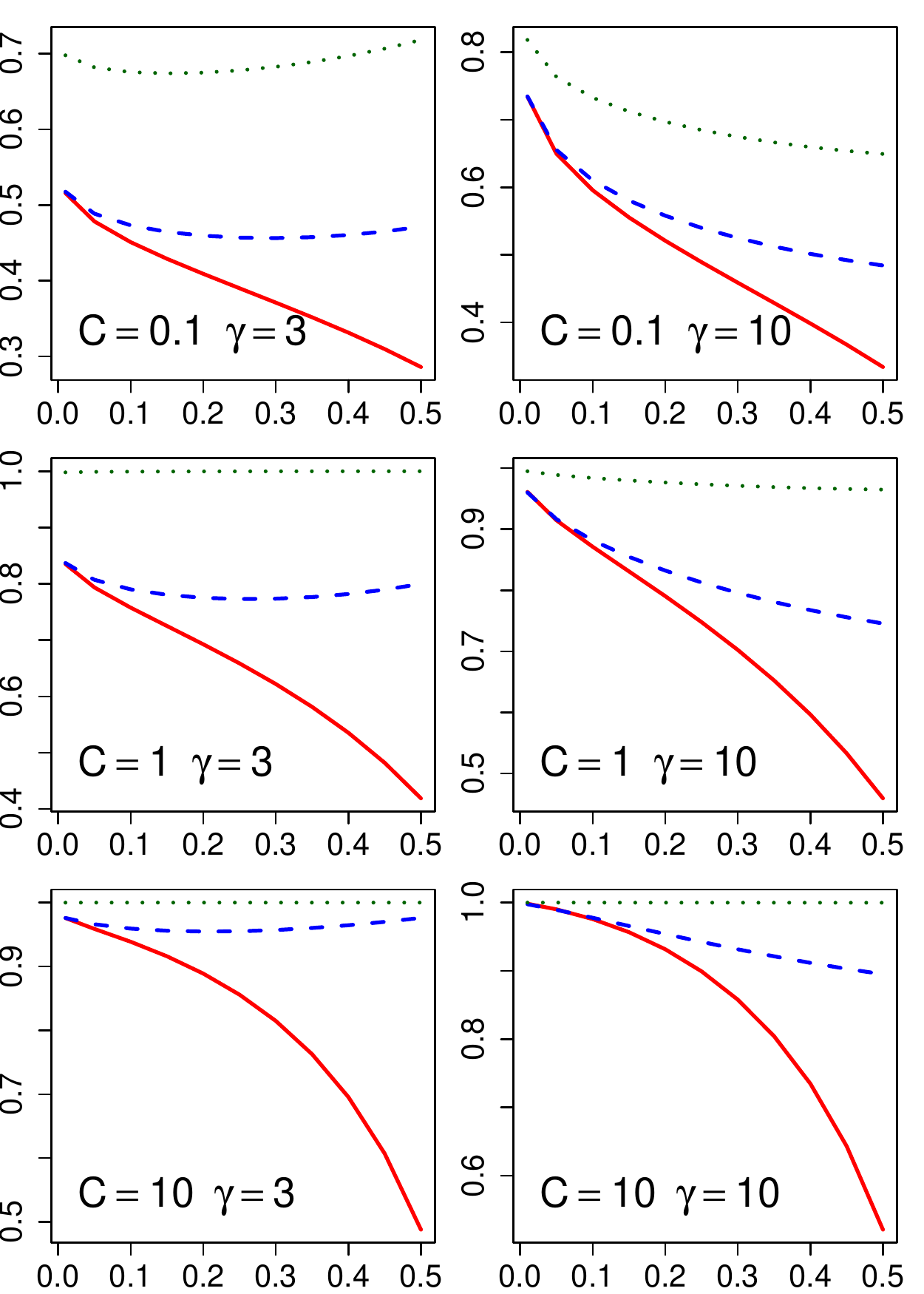}	
	\end{subfigure}
\caption{Plots of MP, P1, and P2 against sparsity parameter $p$ for Pareto distributions. The red solid, blue dashed, and green dotted lines represent the Bayes oracle, the $\alpha_\infty$-BH procedure and the $1/\log(m)$-BH procedure respectively. In panel (a), the diagonal lines are the MP=$p$ line. The ranges of the vertical axes are different across the plots of P1 and P2.}\label{fig:pareto2}
\end{figure}

 \smallskip
The average number of misclassified observations, type \rom{1}, and type \rom{2} errors based on 1000 replicates were used to estimate the misclassification probability (MP), probability of type \rom{1} errors (P1), and probability of type \rom{2} errors (P2), respectively. Figures \ref{fig:t2} and \ref{fig:pareto2} display the plots of the three error measurements against $p$ for Student's t and Pareto distributions, respectively. The solid, dashed, and dotted lines, respectively, represent the Bayes oracle, the $\alpha_\infty$-BH procedure, and the $1/\log(m)$-BH procedure. Here, the $\alpha_\infty$-BH procedure and the Bayes oracle behave similarly if $p$ is small. The solid lines and the dashed lines are almost identical in most situations when $p < 0.1$. When $p$ is large, the $\alpha_\infty$-BH procedure has lower P1 and higher P2 than the Bayes oracle, which suggests the former is conservative in identifying signals when they are abundant. Second, in terms of MP, for Student's t distributions, $1/\log(m)$-BH procedure works better when $\gamma$ is larger since Student's t distributions with bigger degrees of freedom are closer to the normal distributions, for which the $1/\log(m)$ level is designed. In general, when applied to polynomial-tailed distributions, the $1/\log(m)$-BH procedure is more conservative in identifying signals than the $\alpha_\infty$-BH procedure. In the plots for $C=1, \gamma=3$, and $C=10$, its P1 is in close vicinity of 0 and P2 is close to 1, which indicates the procedure identifies almost all observations as noises. In the MP panels of both figures, as $C$ grows, the line corresponding to the Bayes oracle lies closer to the $\mbox{MP}=p$ line indicating increasing difficulty in multiple testing problem. This agrees with our findings in the first scenario.

\par \bigskip

\setcounter{chapter}{6}
\setcounter{equation}{0} 
\noindent {\bf 6. Discussion}
\par \smallskip
This paper establishes some asymptotic optimality properties of several multiple testing procedures in a Bayesian framework where the data are generated from distributions with polynomial tails. In particular, it is shown that some of the classical multiple testing procedures attain asymptotically the Bayes oracle property under sparsity. To the authors' knowledge, Theorem 2 is the first result clearly providing an approach to simplify the problem of finding an implementable random ABOS thresholding procedure to the construction of an appropriate estimator of a fixed ABOS threshold. Future work might extend results beyond polynomial-tailed distributions.

In Section 4.1, we show that, for a fixed $m$, the lower bound of the \textsc{BFDR} level that can be controlled is $\beta^*$, see \eqref{eq:beta_our}, with 
\begin{equation}\label{eq:betastar_inf}
\beta^*_\infty = \lim_{m\rightarrow\infty}\beta^* = (1 +  \delta_\infty / C_0)^{-1}.
\end{equation}
In Section 3 of \citet*{chi2007performance}, it is shown that if both $p$ and $u$ do not vary with $m$ and the cdf of the p-value is strictly concave, then as $m$ grows to infinity, the \textsc{BFDR} is always bounded below by 
\begin{equation}\label{eq:beta_chi}
\beta_*= \frac{1-p}{1-p + p\lim_{x\rightarrow \infty}\rho(x)},
\end{equation}
where $\rho(x)$ in our notation is $(1+u)^{-1/2}d(x/\sqrt{1+u})/d(x)$ with limit $(1+u)^{\gamma/2}$. Thus $\beta^*$ in \eqref{eq:beta_our} and $\beta_*$ in \eqref{eq:beta_chi} have the same expression, although they are derived in different contexts. \citet*{chi2007performance} also proved that, under his setting, there is a critical value $\alpha_* > 0$ for the target \textsc{FDR} control level $\alpha$. If $0 < \alpha < \alpha_*$, the power of a multiple testing procedure decays to 0 as $m \rightarrow \infty$ and the \textsc{BFDR} converges to $\beta_*$. Under our setting, we believe that the criticality phenomenon still exists with both $\alpha_*$ and $\beta_*$ replaced by $\beta_\infty^*$ defined in \eqref{eq:betastar_inf}. In panel (c) of Figures \ref{fig:t} and \ref{fig:pareto}, we plot the probability of type \rom{2} error (P2) of the Benjamini-Hochberg procedure against \textsc{FDR} level $\alpha$ when $m=10^6$. The solid vertical lines represent $\beta^*_\infty$. In the plots for $C=1, \gamma=3$, and $C=10$, P2 is very close to 1, which suggests that the power is close to 0.

The asymptotic optimal \textsc{FDR} level of the Benjamini-Hochberg procedure, $\alpha_\infty$, depends on the difficulty index $C$ that is usually unknown. For practical use, according to our simulation results, when $\gamma$ is large, $1/\log(m)$ is a good surrogate for $\alpha_\infty$ in terms  of risk ratios and misclassification probabilities. Although in some situations $1/\log(m)$ is not close to $\alpha_\infty$, the risk ratio does not considerably deviate from 1. It is shown in Figures \ref{fig:t} and \ref{fig:pareto} that the risk ratios are less sensitive to the choice of $\alpha$ for a larger $C$. Therefore, a smaller value, say 0.1, is a safe guess for $C$ when $\gamma$ is small. There could be more delicate methods to estimate $C$. To illustrate an example, let $m_1$ and $m_0$ denote the number of observations with absolute values in intervals $(b, +\infty)$ and $(a_1, a_2)$, respectively. If $b$ is large enough and $a_1, a_2$ are relatively small, by an idea similar to the one used to estimate $\pi_0$ in Chapter 4.5 of \citet*{efron2012large}, we could assume that almost all the $m_1$ observations are signals and almost all the $m_0$ observations are noises. Then 
${m_1}/{m} \approx 2p\{1-D(b/\sigma_1)\},$ and ${m_0}/{m} \approx 2(1-p)\{D(a_2/\sigma_0)-D(a_1/\sigma_0)\}.$
Taking the ratio of these two and using the polynomial tail equivalence of the MPT distribution, we have
\begin{equation*}
\frac{m_0}{m_1} \approx \frac{1-p}{p}\left(\frac{\sigma_0}{\sigma_1}\right)^\gamma\left[\left(\frac{b}{a_1}\right)^\gamma - \left(\frac{b}{a_2}\right)^\gamma\right] \approx \delta^{-1}vu^{-\gamma/2}[(b/a_1)^\gamma-(b/a_2)^\gamma].
\end{equation*}
As $vu^{-\gamma/2} \rightarrow C_0$, $C_0$ could be estimated by 
\begin{equation}
\hat{C}_0 = (\delta_\infty m_0 / m_1)[(b/a_1)^\gamma-(b/a_2)^\gamma]^{-1}.
\end{equation}
Since $C_0$ is an increasing function in $C$, the estimate of $C$ can be solved analytically or numerically depending on the form of $d$. A problem of this method is how to choose $a_1$, $a_2$. We want $a_1$ and $a_2$ to be small enough so that nearly all the observations in intervals $(-a_2, -a_1)$ and $(a_1, a_2)$ are noises. At the same time, $a_1$ and $a_2$ should be large enough so that the polynomial approximation of $D(a_i/\sigma_0), i=1, 2,$ is accurate. In some simulations not shown here, there is no simple solution to this problem.

As far as we know, theories of multiple testing problems for polynomial-tailed distributions have not been well developed in literature. It is unclear whether some of the interesting results for normal distributions still exist for polynomial-tailed distributions. For example, \citet*{bogdan2011asymptotic} mention that their assumption $u \varpropto -\log(p)$ (obtained from ${2\log(v)}/{u} \rightarrow C$ when $\log(\delta) = o(\log(p))$) can be related to asymptotically least-favorable configurations for $l_0[p]$ balls, discussed in \citet*{abramovich2006special}. To examine whether our assumption $vu^{-\gamma/2} \rightarrow C_0$ has similar connection with minimax estimation, a vital question to be answered is what the configurations for polynomial-tailed distributions look like.

Global-local shrinkage priors have received much attention recently in Bayesian analysis. \citet*{ghosh2013asymptotic} showed that a multiple testing procedure based on a group of global-local shrinkage priors can asymptotically achieve the oracle Bayes risk up to a multiplicative constant. In the same vein, we would like to examine, in future work, whether and how global-local shrinkage priors can be used for polynomial-tailed distributions.
\par \bigskip

\setcounter{chapter}{7}
\setcounter{equation}{0} 
\noindent {\bf 7. Supplementary Material}
\par \smallskip
The online supplementary material includes proofs of our results.
\par \bigskip

\noindent {\large\bf Acknowledgment}
\par \smallskip
Ghosh's research was supported in part by the NSF Grant SES-1327359. Li's research is supported by the Fundamental Research Funds for the Central Universities JBK120509 and JBK140507.
\par \bigskip

\bibliographystyle{sinica}
\bibliography{ABOS_poly}

\vskip .65cm
\noindent
Department of Statistics, University of Florida, Gainesville, Florida 32611-8545, U.S.A.
\vskip 2pt
\noindent
E-mail: xytang@stat.ufl.edu
\vskip 2pt
\noindent
School of Statistics \& Center of Statistical Research, Southwestern University of Finance and Economics, Chengdu, Sichuan, China
\vskip 2pt
\noindent
E-mail: likec@swufe.edu.cn
\vskip 2pt
\noindent
Department of Statistics, University of Florida, Gainesville, Florida 32611-8545, U.S.A.
\vskip 2pt
\noindent
E-mail: ghoshm@stat.ufl.edu
\vskip .3cm
\end{document}